\documentclass[11pt]{amsart}
\usepackage{breqn,color}
\usepackage{amsmath,amsthm,amssymb,amsfonts,mathtools,bm,bbm,array,float,mathdots,dsfont,mathrsfs,latexsym}
\usepackage{caption,subcaption,multirow,longtable,lscape,wrapfig,tikz,tikz-cd,comment,pgffor}
\usepackage[colorlinks=true, pdfstartview=FitV,linkcolor=blue,citecolor=blue,urlcolor=blue]{hyperref}
\usepackage{fullpage}
\usepackage{xcolor}
\usepackage[colorinlistoftodos]{todonotes}
\usepackage{verbatim}
\usepackage{colonequals}

\usepackage{numprint}
\npthousandsep{,}  

\usepackage{diagbox}
\usepackage{booktabs} 
\usepackage{multirow}
\allowdisplaybreaks

\usepackage[square,sort,comma,numbers]{natbib}
\defcitealias{lmfdb}{LMFDB}
\iftrue

\fi

\usepackage{afterpage}

\DeclareMathOperator{\ord}{ord}

\theoremstyle{definition}

\def\Z{\mathbb{Z}}

\usepackage[normalem]{ulem}

\newenvironment{red}{\relax\color{red}}{\hspace*{.5ex}\relax}
\newenvironment{blue}{\relax\color{blue}}{\hspace*{.5ex}\relax}
\newcommand{\ber}{\begin{red}}
\newcommand{\er}{\end{red}}
\newcommand{\beb}{\begin{blue}}
\newcommand{\eb}{\end{blue}}
\newcommand{\seteq}{\mathbin{:=}}

\numberwithin{equation}{section}
\numberwithin{figure}{section}
\numberwithin{table}{section}
\allowdisplaybreaks

\begin{document}

\title{Machine learning the vanishing order of \\ rational $L$-functions}

\author[J. Bieri]{Joanna Bieri}
\address{University of Redlands, Redlands, CA 92373, USA}
\email{\href{mailto:joanna_bieri@redlands.edu}{joanna\_bieri@redlands.edu}}

\email{}\author[G. Butbaia]{Giorgi Butbaia}
\address{Department of Physics and Astronomy, University of New Hampshire, Durham, NH 03824,
USA}
\email{\href{mailto:giorgi.butbaia@unh.edu}{giorgi.butbaia@unh.edu}}

\author[E. Costa]{Edgar Costa}
\address{
  Department of Mathematics,
  Massachusetts Institute of Technology,
  Cambridge,
  MA 02139,
  USA
}
\email{\href{mailto:edgarc@mit.edu}{edgarc@mit.edu}}
\urladdr{\url{https://edgarcosta.org}}

\author[A. Deines]{Alyson Deines}
\address{Center for Communications Research, La Jolla, USA}
\email{\href{mailto:aly.deines@gmail.com}{aly.deines@gmail.com}}

\author[K.-H. Lee]{Kyu-Hwan Lee}
\address{Department of Mathematics, University of Connecticut, Storrs, CT 06269, USA  \hfill \break \indent Korea Institute for Advanced Study, Seoul 02455, Republic of Korea}
\email{\href{mailto:khlee@math.uconn.edu}{khlee@math.uconn.edu}}

\author[D. Lowry-Duda]{David Lowry-Duda}
\address{ICERM, Providence, RI, 02903, USA}
\email{\href{mailto:david@lowryduda.com}{david@lowryduda.com}}
\urladdr{\url{https://davidlowryduda.com}}

\author[T. Oliver]{Thomas Oliver}
\address{University of Westminster, London, UK}
\email{\href{mailto:T.Oliver@westminster.ac.uk}{T.Oliver@westminster.ac.uk}}

\author[Y. Qi]{Yidi Qi}
\address{Department of Physics, Northeastern University, Boston, MA, USA \hfill \break \indent NSF Institute for Artificial Intelligence and Fundamental Interactions, Cambridge, MA, USA}
\email{\href{mailto:y.qi@northeastern.edu}{y.qi@northeastern.edu}}

\author[T. Veenstra]{Tamara Veenstra}
\address{Center for Communications Research, La Jolla, USA}
\email{\href{mailto:tamarabveenstra@gmail.com}{tamarabveenstra@gmail.com}}

\date{\today}

\begin{abstract}
In this paper, we study the vanishing order of rational $L$-functions from a data scientific perspective.
Each $L$-function is represented in our data by finitely many Dirichlet coefficients, the normalisation of which depends on the context. 
We observe murmuration-like patterns in averages across our dataset,
find that PCA clusters rational $L$-functions by their vanishing order, and record that LDA and neural networks may accurately predict this quantity. 
\end{abstract}

\maketitle

\section{Introduction}\label{s.intro}

In this paper we study rational $L$-functions.
Notable examples include those arising from Artin representations, 
arithmetic curves (such as elliptic curves and genus $2$ curves defined over number fields), and certain modular forms (such as classical, Hilbert, and Bianchi modular forms).
We apply techniques from machine learning to the Dirichlet coefficients of each of these $L$-functions.

We will use techniques such as principal component analysis, feed-forward neural networks, and convolutional neural networks.
We are mainly interested in the vanishing order at the central point, which is conjecturally related to the underlying arithmetic, as per the conjectures of Birch and Swinnerton-Dyer and, more generally, Beilinson--Bloch--Kato. 

Following this introduction, Section~\ref{s:def} provides the necessary definitions and examines various patterns in the averages of the coefficients. These are closely related to murmurations.
Murmurations were first observed in the context of elliptic curves over $\mathbb{Q}$~\cite{HLOP}.
Section~\ref{s:pca} applies principal component analysis, linear discriminant analysis, and neural networks to datasets of rational $L$-functions, aiming to learn their vanishing orders. 
Comparable techniques were applied to arithmetic curves in \cite{HLOa,HLOc,HLOP}, to number fields in \cite{HLOb,AHLOS}, and to a tiny heterogeneous set of rational $L$-functions in \cite{O24}.
While neural networks have been applied to elliptic curve datasets in \cite{KV22} and \cite{Poz}, the novelty of this work lies in extending these techniques to a large, heterogeneous dataset comprising $L$-functions with diverse origins. 
Furthermore, we explore transfer learning from one sub-dataset to another.
In Section~\ref{s.conc} we summarise our main results and offer some future directions. In Section~\ref{sect:appendix}, we provide an appendix of murmuration patterns for data outside of the main dataset studied in the main text.

\subsection*{Acknowledgements}
Costa was supported by Simons Foundation grants 550033 and SFI-MPS-Infrastructure-00008651.
Lee was supported by Simons Foundation grant 712100.
Lowry-Duda was supported by Simons Foundation grant 546235.
Qi was supported by the NSF grant PHY-2019786 (the NSF AI Institute for Artificial Intelligence and Fundamental Interactions).
We would also like to express our gratitude to the organizers of the Mathematics and Machine Learning Program at CMSA during the Fall 2024 semester, where this project was initiated.

\section{Preliminaries}\label{s:def}
In this section, we define the terminology used, describe our primary dataset, and document various murmuration-like patterns in averages across the dataset.

\subsection{Definitions}\label{s.def}

Loosely speaking, an $L$-function is a Dirichlet series
$L(s) \seteq \sum_{n=1}^{\infty} a_n n^{-s}$  that has an Euler product and satisfies a functional equation.
Each $L$-function is a generating function for the sequence of coefficients
$\{ a_n \}_{n=1}^\infty$ and the analytic properties of $L(s)$ relate to
arithmetic properties of the coefficients.
A precise axiomatisation was suggested by Selberg~\cite{S92}.

In greater detail, the functional equation for an $L$-function is easier to describe after multiplying by a function of the form
$N^s\prod_{i=1}^k\Gamma(\omega_is+\mu_i)$, in which $N\in\mathbb{Z}$, $\omega_i\in\mathbb{R}_{>0}$ and $\mu_i\in\mathbb{C}$ with $\mathrm{Re}(\mu_i)>0$.
We refer to the integer $N$ as the \emph{conductor} of $L(s)$ and $L_\infty(s) = \prod_{i = 1}^k \Gamma(\omega_i s + \mu_i)$ as the factor \emph{at infinity}.
The product $\Lambda(s) = N^s L_\infty(s) L(s)$ is called the completed $L$-function and, for each $L$-function we consider in this paper, satisfies a functional equation of the shape $\Lambda(s) = \varepsilon \Lambda(\delta - s)$ for some $\delta > 0$ and $\varepsilon$ with $\lvert \varepsilon \rvert = 1$.

The degree of $L(s)$ is given by the expression $d=2\sum_{i=1}^k\omega_i$ and is conjectured to be an integer for $L$-functions arising naturally in number theory.
(In this paper, we consider $L$-functions of degrees $1$, $2$, and $4$).
The product of two $L$-functions is another $L$-function.
An $L$-function is said to be \emph{primitive} if it cannot be written as a product of $L$-functions with smaller degree. 
To each $L$-function we associate a measure of complexity called the \emph{analytic conductor}, given by
\[
  A
  =
  N\cdot\exp \left(
    2\mathrm{Re}\left(\frac{L'_{\infty}(1/2)}{L_{\infty}(1/2)}\right)
  \right),
\]
where $L_{\infty}(s)$ is the factor at infinity as defined above.
To compare across different degrees, we also use the \emph{root analytic conductor}, which is the $d$th root of the analytic conductor if $L$ has degree $d$.
(With this, $L(s)$ and $L^2(s)$ have the same root analytic conductor).

We refer to the sequence $\{ a_n \}_{n=1}^{\infty}$ as the Dirichlet coefficients of $L(s)$.
The Dirichlet coefficients of an $L$-function are multiplicative in the sense that, if $\mathrm{gcd}(m,n)=1$, then $a_ma_n=a_{mn}$.
Though we refer to the coefficients as $a_n$, there is a choice of normalisation.
For example, if $L(s)$ satisfies a functional equation $\Lambda(s) = \Lambda(\delta - s)$, then multiplying $a_n$ by $\sqrt{n}$ gives another $L$-function with a functional equation of the shape $\Lambda(s) = \Lambda\bigl((1 + \delta) - s\bigr)$.
After normalising through multiplication by a certain power $n^\alpha$, we can put an $L$-function in the \emph{analytic normalisation} $L_{\textup{an}}(s)$, wherein the functional equation has shape $s \mapsto 1 - s$.

We say that $L(s)$ is an \emph{arithmetic} $L$-function if its Dirichlet coefficients are algebraic numbers.
For an arithmetic $L$-function $L(s)$ given in its analytic normalisation $L_{\textup{an}}(s)$, we define its \emph{motivic weight} to be the minimal $w \in \mathbb{Z}_{\geq 0}$ such that, for all $n \geq 1$, the number $a_n n^{w / 2}$ is an algebraic integer.
In such a case we also define the \emph{coefficient field} of $L(s)$ to be the
number field generated by $\{ a_n n^{w/2} \}_{n = 1}^\infty$.
If the coefficient field is $\mathbb{Q}$, then we refer to $L(s)$ as a \emph{rational} $L$-function.
For example, Dirichlet characters (resp. elliptic curves over $\mathbb{Q}$) determine rational $L$-functions with motivic weight $0$ (resp. $1$).

Different normalisations can make certain number theoretic properties clearer.
Taking an arithmetic $L$-function in the analytic normalisation
$L_{\textup{an}}(s)$ and scaling by the motivic weight by multiplying each $a_n$ by
$n^{w/2}$ gives the \emph{arithmetic normalisation}, which we refer to as $L(s)$
(so that $L_{\textup{an}}(s) = L(s + \frac{w}{2})$).
In this normalisation, the coefficients are algebraic integers, the functional
equation has the shape $s \mapsto 1 + w - s$, and the central point is
$s = (w + 1)/2$.
If $L(s)$ is rational, then, with the arithmetic normalisation, its Euler product can be written as a product over rational primes
$
L(s)=\prod_pL_p(p^{-s})^{-1},
$
with $L_p(T)\in \Z[T]$. 

As mentioned in Section~\ref{s.intro}, the vanishing order of an $L$-function at its central point $r \colonequals \ord_{s  = (w + 1)/2} L(s)$ is expected to reflect underlying arithmetic structure.
For example, the conjecture of Birch and Swinnerton-Dyer asserts that the vanishing order of the $L$-function attached to an abelian variety defined over a number field is equal to the rank of that abelian variety as a finitely generated abelian group.
Moreover, in the arithmetic normalisation, the leading Taylor coefficient is conjecturally given by a formula involving other arithmetic quantities.

\subsection{Describing the dataset}

\begin{figure}[ht!]
\centering
\includegraphics[width=0.9\textwidth]{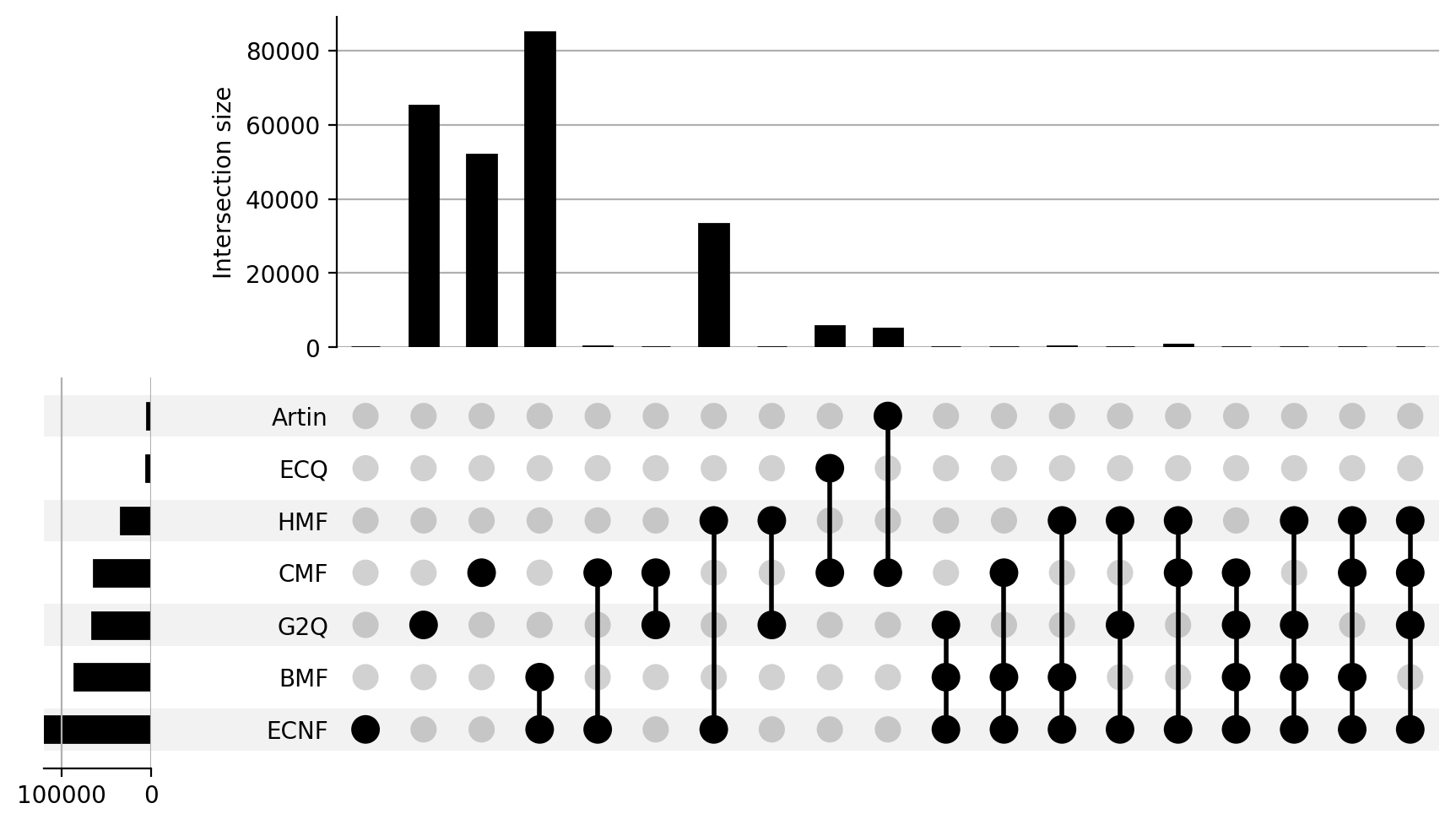}
\caption{UpSet plot for \texttt{RAT}. Each row corresponds to the named subset, and each column corresponds to an intersection of subsets. A black circle is used to indicate that a subset is involved in a particular intersection (the vertical lines simply indicate the orientation of the plot). The horizontal bars show the number of datapoints in each subset (row), and a vertical bar above shows the size of an exclusive intersection, i.e., it represents the number of elements unique to that specific intersection.} 
\label{fig:upset}
\end{figure}

Throughout this paper, we use a dataset \texttt{RAT}~\cite{rat} consisting of rational $L$-functions taken from \cite{lmfdb}, and in particular based on work from~\cite{cremona1981modular, cremona1984hyperbolic, whitley1990modular, cremona1997algorithms, boberdatabase, jones2017artin, booker2016database, donnelly2021database, best2021computing, cremona2021q}.
The \texttt{RAT} dataset contains information about
\numprint{248359} rational $L$-functions with root analytic conductor less than 4.
We restrict to root analytic conductor $< 4$ because this dataset contains a relatively balanced set of $L$-functions from different origins in this range.
The available data from the LMFDB becomes more skewed for larger root analytic conductor.
The \texttt{RAT} dataset includes the sub-datasets \texttt{CMF} (classical modular forms), \texttt{ECNF} (elliptic curves defined over number fields), \texttt{DIR} (Dirichlet characters) and \texttt{G2Q} (genus 2 curves defined over $\mathbb{Q}$).
Within \texttt{CMF}, we find familiar sub-datasets such as \texttt{ECQ} (elliptic curves over $\mathbb{Q}$) and \texttt{ART} (Artin representations\footnote{There is one exception. The Riemann zeta function is contained in \texttt{ART} but not in \texttt{CMF}.}).
Similarly, within \texttt{ECNF} we find \texttt{BMF} (Bianchi modular forms) and \texttt{HMF} (Hilbert modular forms).
In Figure~\ref{fig:upset}, we present an UpSet plot for \texttt{RAT}, which is a method for visualising data with more than three intersecting subsets (cf.~\cite{lex2014upset}).

Contained within the \texttt{RAT} dataset is a sub-dataset \texttt{PRAT} consisting of \numprint{186114} primitive rational $L$-functions.
Since imprimitive $L$-functions can be written as a product of primitive $L$-functions, we expect \texttt{RAT} to contain theoretically redundant information (for example, all the classical modular forms that do not arise from elliptic curves). 
We caution the reader that this may not be entirely true in practice, since there is some chance that the imprimitive $L$-functions in $\texttt{RAT}$ do not factor into primitive $L$-functions that are included in $\texttt{PRAT}$. 
Nevertheless, for our experiments, we will restrict to \texttt{PRAT}, in order to exclude the possibility of redundant information.
There are other idiosyncrasies in the \texttt{PRAT} dataset.
For example, there is no overlap between \texttt{CMF} and \texttt{ECNF} (although this is not true in general). The number of datapoints in each sub-dataset is given in Table~\ref{tab:counts}, and the number of datapoints in various pairwise intersections is given in Table~\ref{tab:intersections}.

In the \texttt{(P)RAT} dataset, each row contains information about one $L$-function and includes several columns~\cite{rat}.
This includes the following properties of the $L$-function:
\begin{itemize}
\item primitivity of the $L$-function;
\item $N$, the conductor of the $L$-function;
\item $w$, the motivic weight of the $L$-function;
\item $d$, the degree of the $L$-function;
\item $\{a_p\}_{p \leq 1000}$, the Dirichlet coefficients (in the arithmetic normalisation) at primes less than $1000$;
\item $r$, the vanishing order at the central point;
\item $t \subseteq \{
    \texttt{CMF} ,
    \texttt{ECQ} ,
    \texttt{ART} ,
    \texttt{ECNF},
    \texttt{BMF} ,
    \texttt{HMF} ,
    \texttt{DIR} ,
    \texttt{G2Q}
\}
$, identifying the type of objects present in \cite{lmfdb} that give rise to this $L$-function (at time of creation of this dataset).
\end{itemize}

In this paper, we will be particularly interested in learning the order of vanishing $r$  column from the Dirichlet coefficients $\{a_p\}_{p \leq 1000}$ .

\begin{table}
\parbox{0.47\linewidth}{
  \begin{center}
    \begin{tabular}{c|c}
    sub-dataset & number of datapoints\\
    \hline
    \texttt{CMF}  & \numprint{9675}  \\
    \texttt{ECQ}  & \numprint{5860}  \\  
    \texttt{ART}  & \numprint{2557} \\ 
    \texttt{ECNF} & \numprint{113489} \\ 
    \texttt{BMF}  & \numprint{81805} \\  
    \texttt{HMF}  & \numprint{31964} \\  
    \texttt{DIR}  & \numprint{274} \\  
    \texttt{G2Q}  & \numprint{62789}\\
    \end{tabular}
    \captionsetup{width=0.85\linewidth}
    \caption{Number of datapoints for various sub-datasets of \texttt{PRAT}.}
    \label{tab:counts}
  \end{center}
}

\hfill
\parbox{0.48\linewidth}{
  \begin{center}
    \begin{tabular}{c|c}
    intersection & number of datapoints\\
    \hline
    \texttt{ART}$\cap$\texttt{CMF}  & \numprint{2556}  \\ 
    \texttt{BMF}$\cap$\texttt{ECNF} & \numprint{81805}  \\ 
    \texttt{BMF}$\cap$\texttt{G2Q}  & \numprint{72} \\  
    \texttt{BMF}$\cap$\texttt{HMF}  & \numprint{272} \\  
    \texttt{CMF}$\cap$\texttt{ECQ}  & \numprint{5860} \\ 
    \texttt{ECNF}$\cap$\texttt{G2Q} & \numprint{113} \\ 
    \texttt{ECNF}$\cap$\texttt{HMF} & \numprint{31956} \\ 
    \texttt{G2Q}$\cap$\texttt{HMF}  & \numprint{50} \\
    \end{tabular}
    \captionsetup{width=0.85\linewidth}
    \caption{Intersection size for various sub-datasets of \texttt{PRAT}.}
    \label{tab:intersections}
  \end{center}
}
\end{table}

\subsection{Murmurations in \texttt{PRAT}} \label{s:murm}
In \texttt{(P)RAT}, the Dirichlet coefficients $\{ a_p \}_{p \leq 1000}$ are given in the arithmetic normalisation.
In the context of murmurations, it is common to (additionally) normalize the Dirichlet coefficients of prime index by defining
\begin{equation}\label{eq.aptilde}
\widetilde{a}_p=\frac{a_p}{{p^{(w-1)/2}}}\in[-d\sqrt{p},d\sqrt{p}].
\end{equation}
We note that the normalisation in equation~\eqref{eq.aptilde} is unnatural to both number theorists and data scientists, but, for historical and aesthetic reasons, it seems to be the right choice for the purposes of this section. 
For example, equation~\eqref{eq.aptilde} specialises to the normalisation in the initial papers and correspondence on murmuration phenomena~\cite{HLOP, drewLetter}.
In this section, we will plot murmuration-like averages for sub-datasets in \texttt{PRAT}.
True murmuration behavior relates averages of coefficients in particular conductor ranges~\cite{sarnakLetter}.
Since our datasets include large ranges of conductor, these figures are not exactly murmurations.
For ease of exposition, we will refer to these plots simply as ``murmurations''.

\begin{figure}[ht!]
\centering
\includegraphics[width=0.9\textwidth]{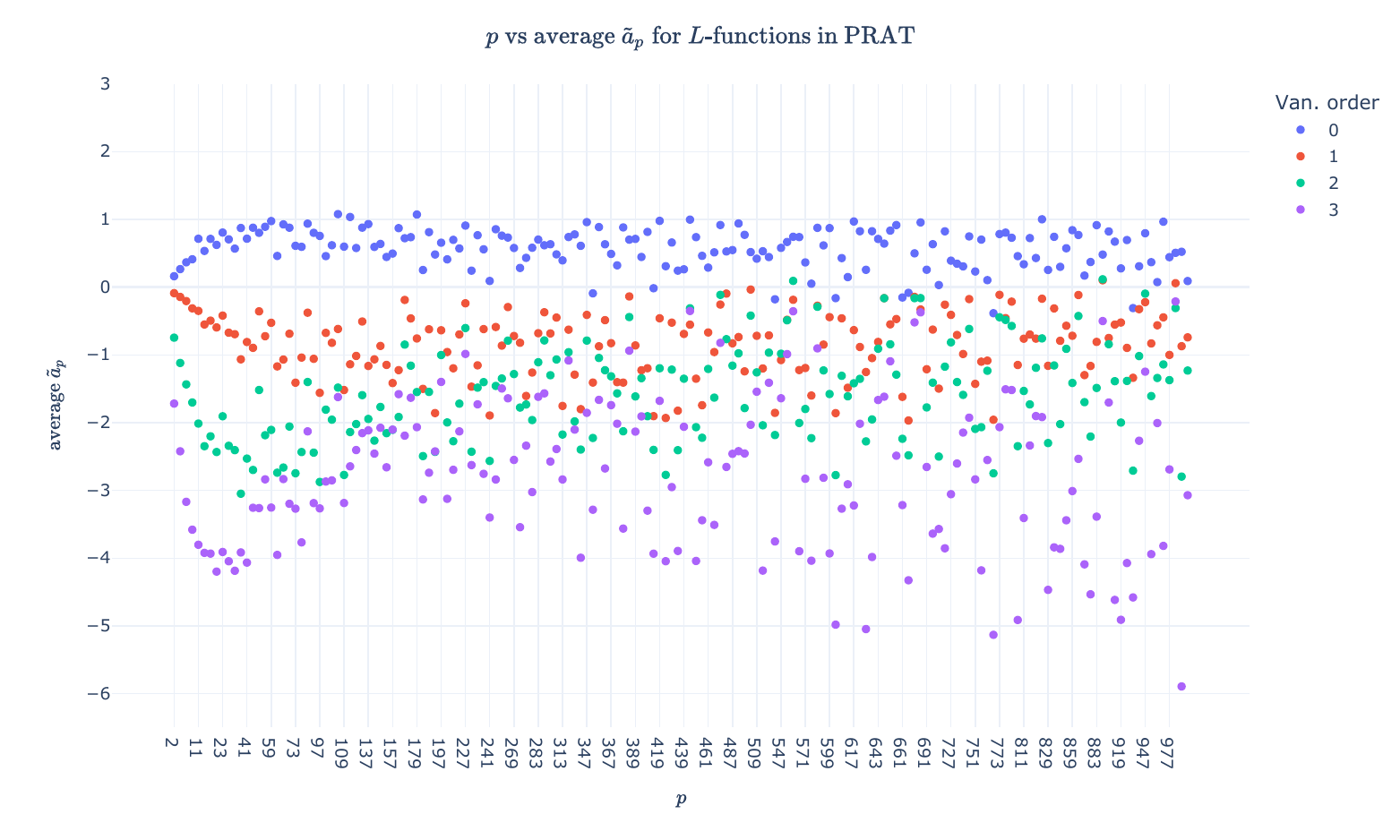}
\caption{Average value of $\widetilde{a}_p$ for primitive rational $L$-functions with specified vanishing order, excluding the $9$ $L$-functions with vanishing order $4$.}
\label{fig:rat-avg-ap-all}
\end{figure}

The counts for each vanishing order in \texttt{PRAT} and varying sub-datasets are recorded in 
Table~\ref{tab:sub-prat-order-counts}.
Since \texttt{PRAT} only contains $9$ data-points with vanishing order equal to $4$, we restrict to the dataset $\texttt{PRAT}_{\leq3}$ consisting of $L$-functions in \texttt{PRAT} with vanishing order $\leq3$.
In Figure~\ref{fig:rat-avg-ap-all},
we plot the average value of $\widetilde{a}_p$
for primitive rational $L$-functions in $\texttt{PRAT}_{\leq3}$ with each vanishing order.

\begin{table}
\begin{tabular}{c|ccccc}
&\multicolumn{5}{c}{vanishing orders}\\
\hline
sub-dataset & $0$ & $1$ & $2$ & $3$ & $4$ \\
\hline
\texttt{PRAT} & \numprint{59863} & \numprint{93657} & \numprint{29748} & \numprint{2837} & $9$\\
\hline
\hline
\texttt{ECNF} & \numprint{42558} & \numprint{61243} & \numprint{9661} & $27$ &\\
\texttt{BMF}  & \numprint{28280} & \numprint{44773} & \numprint{8724} & $26$ & \\
\texttt{HMF}  & \numprint{14443} & \numprint{16582} & \numprint{938}  & $1$ &\\
\texttt{G2Q}  & \numprint{10827} & \numprint{29155} & \numprint{19988}& \numprint{2810} &$9$\\
\texttt{CMF}  & \numprint{6245} & \numprint{3330}  & \numprint{100}  && \\
\texttt{ECQ}  & \numprint{2901} & \numprint{2862}  & \numprint{97}   &&\\
\texttt{Artin}& \numprint{2557} &&&& \\
\texttt{DIR}  & \numprint{274} &&&& \\
\end{tabular}

\caption{Counts for different orders of vanishing in sub-datasets of \texttt{PRAT}.}

\label{tab:sub-prat-order-counts}
\end{table}

In Table~\ref{tab:prat-mw-counts} (resp.~\ref{tab:prat-d-counts}), we record the number of datapoints in $\texttt{PRAT}_{\leq3}$ with given motivic weight (resp. degree).
We observe that the vast majority of $L$-functions in this dataset have degree $4$ and motivic weight $1$.
Restricting to datapoints in $\texttt{PRAT}_{\leq3}$ with this degree and motivic weight excludes a relatively small number of outliers and leaves behind a dataset $\texttt{PRAT}^{\star}$ containing \numprint{176156} $L$-functions, all of which are instances of either \texttt{ECNF} or \texttt{G2Q} (and some of which are both).   
In particular, the dataset \texttt{ECQ} is not contained within $\texttt{PRAT}^{\star}$.
On the other hand, we note that \texttt{ECQ} is already extensively studied in the literature \cite{HLOc,HLOP,Poz,KV22}.
(For comparison, we examine the $L$
-functions in \texttt{PRAT} but not in $\texttt{PRAT}^{\star}$
in Section~\ref{sect:appendix}.)

\begin{table}
\parbox{0.48\linewidth}{\begin{center}
\begin{tabular}{c|c}
 $w$ & count \\
 \hline
$1$   & \numprint{182025}\\
$0$   &   \numprint{2830}\\
$3$   &    \numprint{544}\\
$2$   &    \numprint{177}\\
$5$   &    \numprint{169}\\
\end{tabular}
\captionsetup{width=0.75\linewidth}
\caption{Counts for different motivic weights in $\texttt{PRAT}_{\leq3}$.}
\label{tab:prat-mw-counts}\end{center}
}\hfill\parbox{0.48\linewidth}{%
\begin{center}\begin{tabular}{c|c}
$d$ & count\\
\hline
$4$   & \numprint{176165}\\
$2$    &  \numprint{9675}\\
$1$   &   \numprint{274}\\
\end{tabular}
\captionsetup{width=0.75\linewidth}
\caption{Counts for different degrees in $\texttt{PRAT}_{\leq3}$.}
\label{tab:prat-d-counts}\end{center}}
\end{table}

In what follows, we analyze the murmuration patterns for $\texttt{PRAT}^\star$ and various sub-datasets. 
For all murmuration plots in this section, $L$-functions with larger vanishing orders yield smaller average values of $a_p$.
This is consistent with \cite[Figure~8]{HLOP}, and is comparable to the theory of Mestre--Nagao sums.

\begin{figure}[htbp!]
                  \centering
            \includegraphics[width=\textwidth]{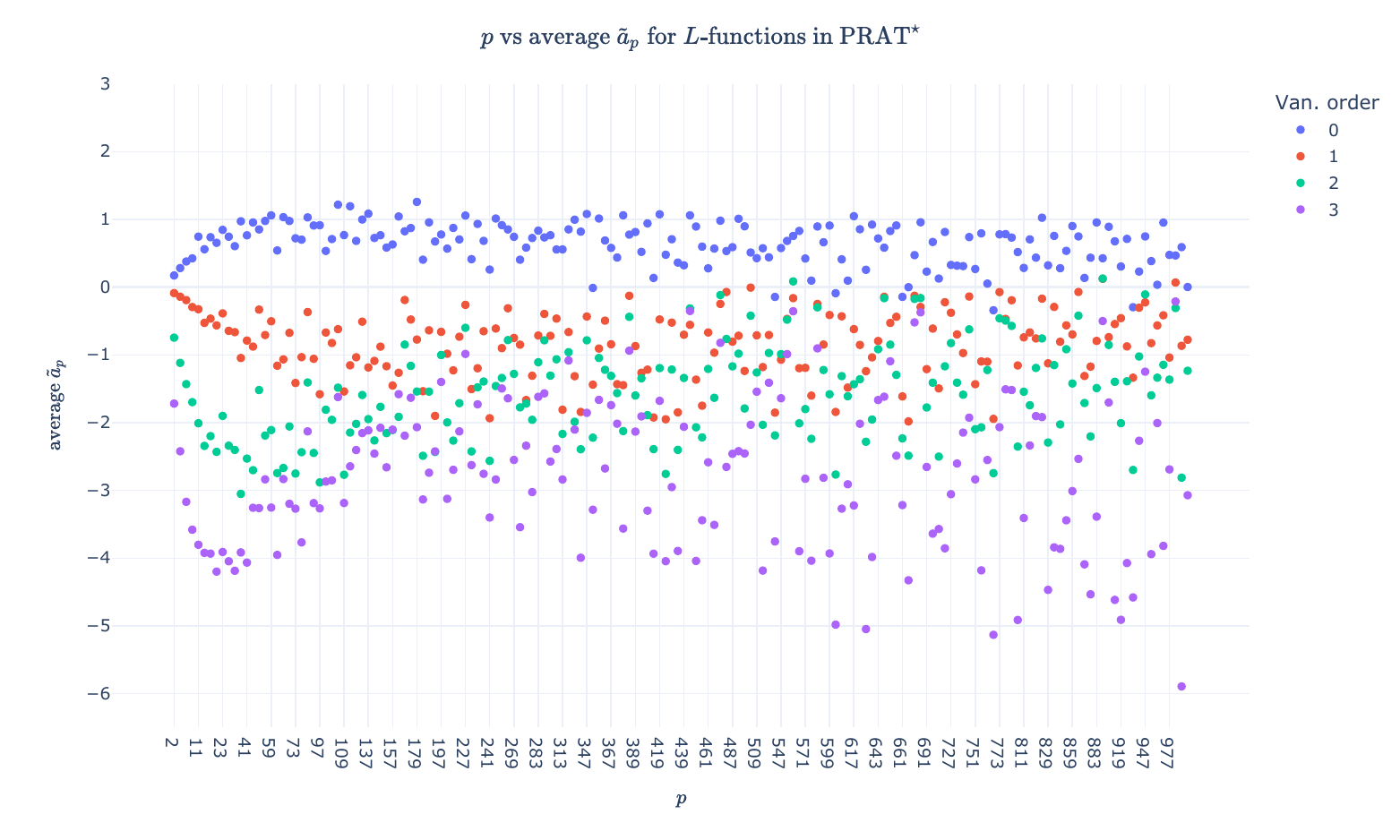}
      
        \caption{Average value of $\widetilde{a}_p$ over $L$-functions with specified vanishing order in
               $\texttt{PRAT}^{\star}$. }
        \label{fig:murm-all}
    \end{figure}

\begin{figure}[htbp!]
        \centering
        \begin{subfigure}[b]{0.48\textwidth}
            \centering
            \includegraphics[width=\textwidth]{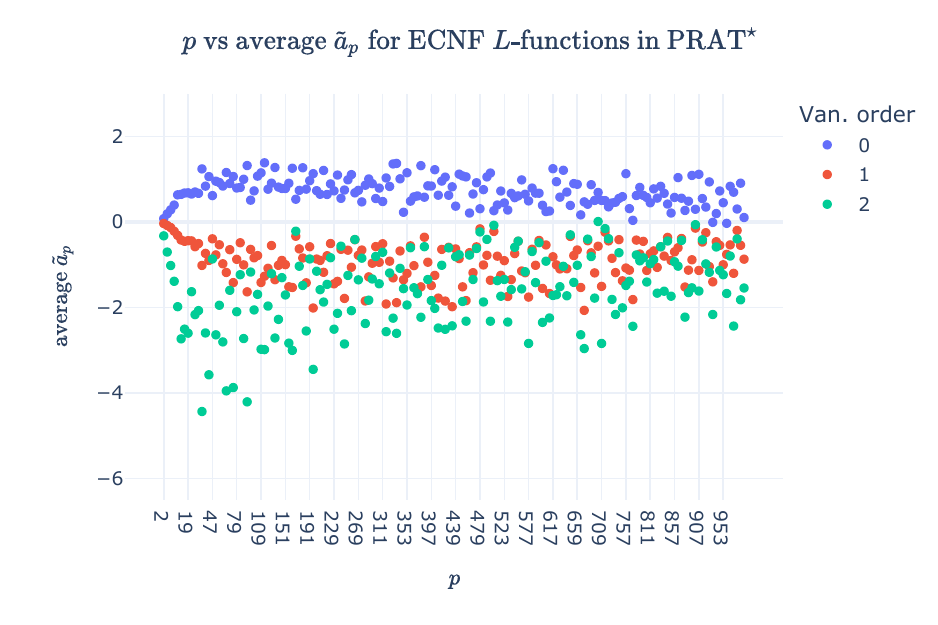}
        \end{subfigure}
        \hfill
        \begin{subfigure}[b]{0.48\textwidth}
            \centering
            \includegraphics[width=\textwidth]{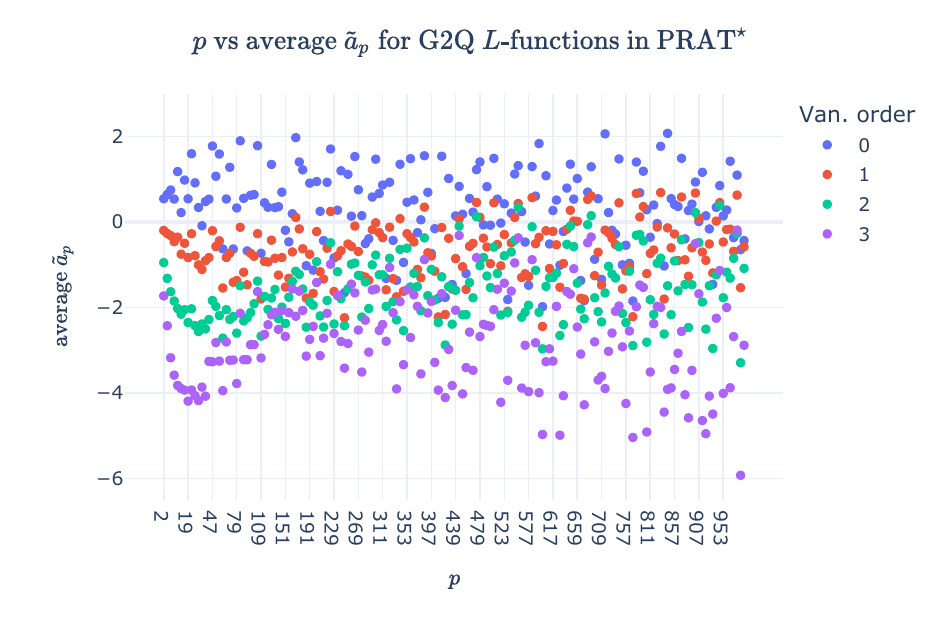}
        \end{subfigure}
        \caption{Average value of $\widetilde{a}_p$ over $L$-functions with specified vanishing order in (left) \texttt{ECNF} and  (right) \texttt{G2Q}
        excluding the $27$ $L$-functions with vanishing order $3$ for \texttt{ECNF}.}
    
        \label{fig:murm-ECNF-vs_G2C}
    \end{figure}

In Figure~\ref{fig:murm-all}, we plot the average value of $\widetilde{a}_p$ for $L$-functions in $\texttt{PRAT}^{\star}$, grouped by vanishing order (since $w=1$, we have $\widetilde a_p=a_p$ in $\texttt{PRAT}^{\star}$).
In Figure~\ref{fig:murm-ECNF-vs_G2C}, we plot \texttt{ECNF} and \texttt{G2Q} separately.
We observe initial separation in the graphs.
More Dirichlet coefficients would allow further exploration of murmuration behavior.

\begin{figure}[ht!]
\centering
\includegraphics[width=0.2\textwidth]{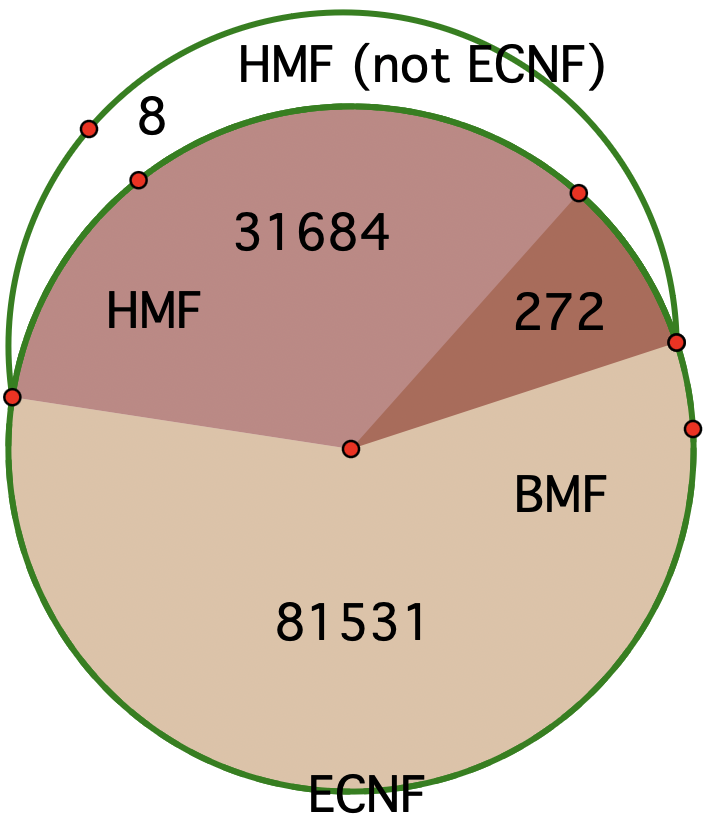}
\caption{Venn diagram for \texttt{ECNF}.}
\label{fig:ECNF-venn-murm}
\end{figure}

\begin{figure}[htbp!]
      \centering
     \begin{subfigure}[b]{0.45\textwidth}
        \centering
       \includegraphics[width=\textwidth]{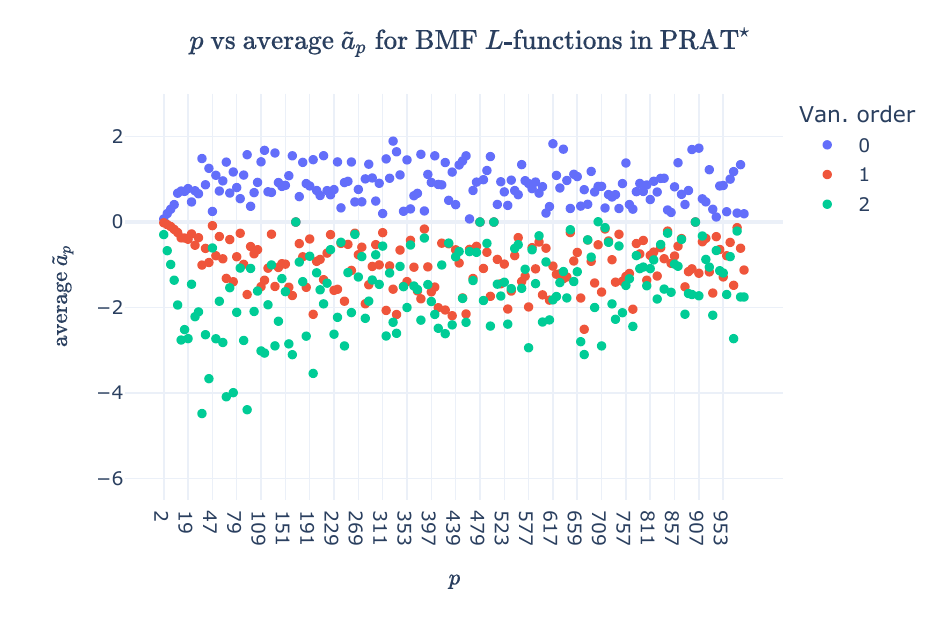}
    \end{subfigure}
        \hfill
        \begin{subfigure}[b]{0.45\textwidth}
            \centering
\includegraphics[width=\textwidth]{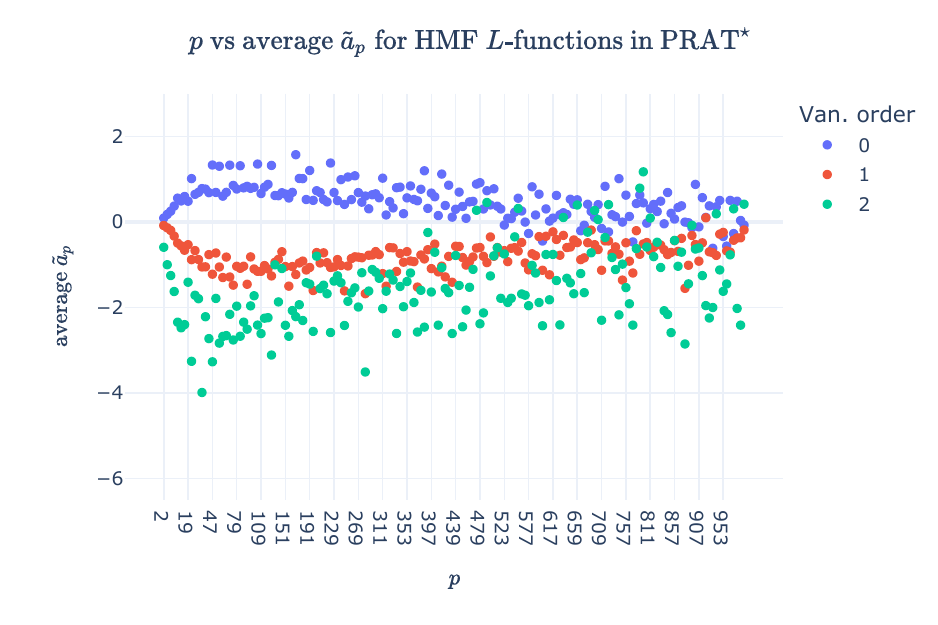}
        \end{subfigure}
        \caption{Average value of $\widetilde{a}_p$ for  (left) \texttt{BMF} and  (right) \texttt{HMF},
        excluding the $26$ (resp. $1$) $L$-functions with vanishing order $3$.} 
        \label{fig:murm-BMF-HMF}
    \end{figure}

In Figure~\ref{fig:ECNF-venn-murm}, we present a Venn diagram for the largest sub-dataset of $\texttt{PRAT}^{\star}$, namely \texttt{ECNF}, which includes \texttt{BMF} and all but $8$ datapoints in \texttt{HMF}. 
In Figure~\ref{fig:murm-BMF-HMF}, we plot the average value of $\widetilde{a}_p$ over \texttt{BMF} and \texttt{HMF}.


\section{Machine learning the vanishing order of rational $L$-functions}\label{s:pca}
We present an unsupervised approach (PCA) and two supervised approaches (LDA and neural networks). 
PCA for elliptic curves was documented in \cite{HLOa,HLOP}, and neural networks were applied to elliptic curves in~\cite{KV22, alexeypozdnyakov2024predicting}.

\subsection{Feature selection and normalisation}
As in Section~\ref{s:murm}, we will define our features using the Dirichlet coefficients column $\{a_p\}_{p \leq 1000}$.
Since the machine learning techniques we will be using are sensitive to the scale of features and we want to prevent 
 features with larger magnitudes from dominating the learning algorithm, we switch normalisation methods.  The normalisation used in Section~\ref{s:murm} produced feature vectors with bounds depending on $d$ and $p$. We now normalize the Dirichlet coefficients with prime index 
by: 
\begin{equation}\label{eq.apbar}
\overline{a}_p = \frac{\widetilde{a}_p}{d \sqrt{p}} = \frac{a_p}{d p^{w/2}}\in[-1,1].
\end{equation}
where $\widetilde{a}_p$ is as in equation~\eqref{eq.aptilde}, $d$ and $w$ are the degree and motivic weight of the $L$-function, see Section~\ref{s.def}.
 Whilst there are several techniques to 
create features with a similar scale, we chose equation~\eqref{eq.apbar} due to its arithmetic nature (reflecting Hasse bounds).
As noted in Section~\ref{s:murm}, the majority of the data has degree 4 and motivic weight 1, so to eliminate outliers we restrict to datapoints in $\texttt{PRAT}^{\star}$. 
For $L$-functions in $\texttt{PRAT}^{\star}$, equation \eqref{eq.apbar} amounts to division by $4\sqrt{p}$.
Each $L$-function in $\texttt{PRAT}^{\star}$ is then represented by the following vector:
\begin{equation}\label{eq.vectors}
v(L)=(\bar{a}_2,\bar{a}_3,\bar{a}_5,\cdots,\bar{a}_{997})\in\mathbb{R}^{168},
\end{equation}
which is indexed by primes $<1000$.
Equation~\eqref{eq.vectors} determines a pointcloud:
\[\mathcal{D}=\{v(L):L\in\texttt{PRAT}^{\star}\}\subset\mathbb{R}^{168}.\]

\subsection{Principal component analysis on $\mathcal{D}$} 

\begin{figure}[htbp]
    \centering
    \includegraphics[width=.4\textwidth]{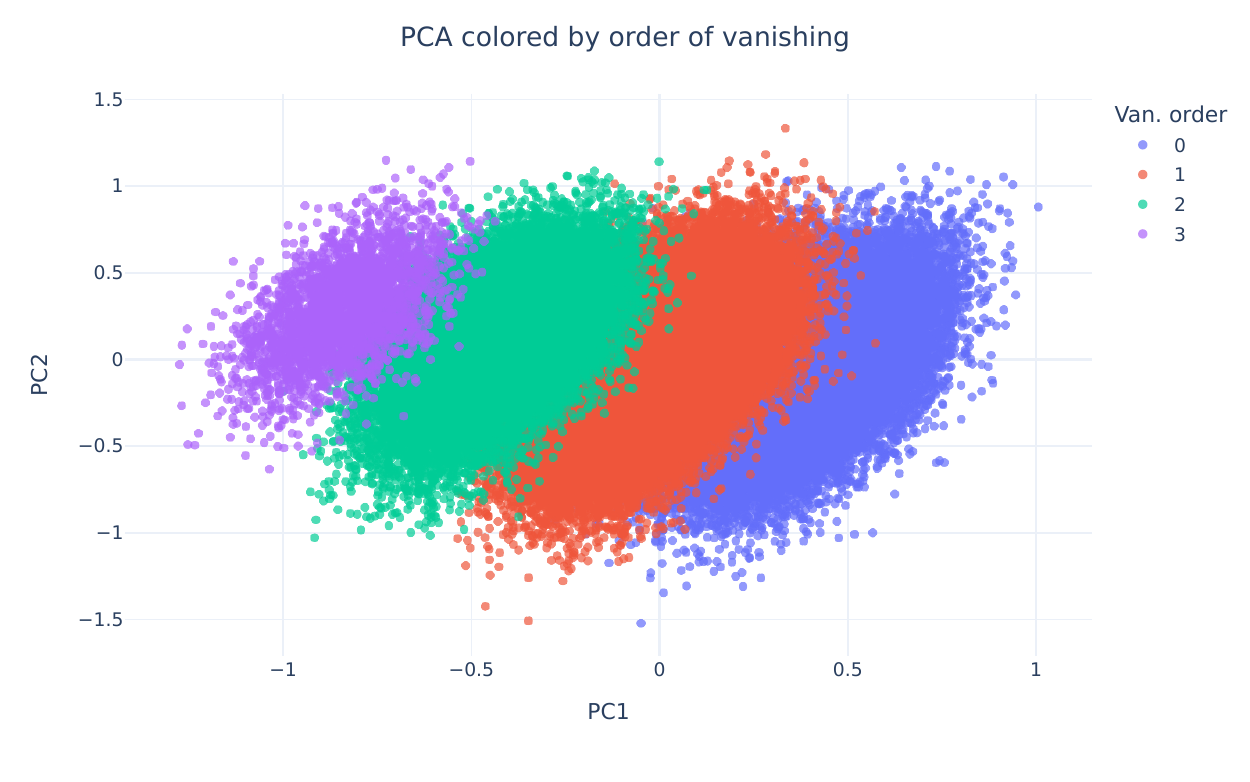}
    \includegraphics[width=.4\textwidth]{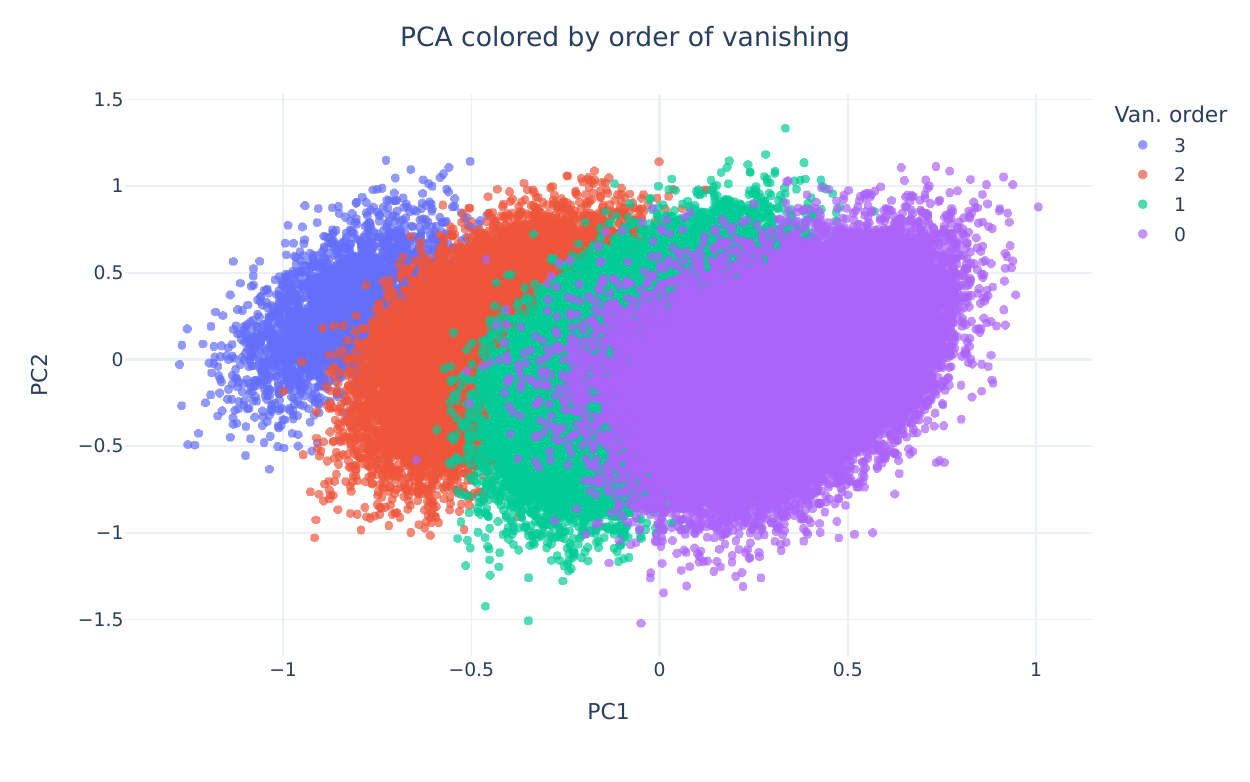}
    \caption{Two-dimensional PCA for the  $\texttt{PRAT}^{\star}$ dataset. Each datapoint is colored by the vanishing order for the underlying $L$-function at its central point.  In the left (resp. right) image we plot the points in ascending (resp. descending) order for the order of vanishing.  }
    \label{fig:2DPCA}
\end{figure}

As depicted in Figure~\ref{fig:2DPCA}, we ran $2$-dimensional PCA on the pointcloud $\mathcal{D}$.
Though some separation is visible, we observe that $2$-dimensional PCA on $\mathcal{D}$ does not fully distinguish the data by vanishing order.
One observes similar behavior for $3$-dimensional PCA. 

The projection onto each principal component is given by a linear combination of $\overline{a_p}$ 
\begin{equation}\label{eq.PCproj}
\sum_{p\leq 997}w_p\overline{a_p}, \ \ w_p\in\mathbb{R}.
\end{equation}
We refer to $w_p$ as the \textsl{weight} at $p$.
We may calculate the weights by diagonalising the covariance matrix for $\mathcal{D}$. 
Equation~\eqref{eq.PCproj} is somewhat analogous to Mestre--Nagao sums of the form 
\begin{equation}\label{eq.mnsum}
S(B)=\frac{1}{\log(B)}\sum_{p<B}\frac{a_p\log(p)}{p},
\end{equation}
which have been used in the detection of high rank  elliptic curves. 
Murmurations for Mestre--Nagao sums were observed in \cite{BKN}.

\subsection{LDA for vanishing order}\label{s.lda}

We used Linear Discriminant Analysis (LDA) to predict the order of vanishing for $L$-functions in subsets of $\texttt{PRAT}^{\star}$ from the feature vectors $v(L)$. 
We divide the data into training and validation sets, with a ratio of $80:20$, stratified by the order of vanishing to ensure that we preserve the proportions of the classes in both the training and validation sets.
Table \ref{LDA-table} shows the results both for the full $\texttt{PRAT}^{\star}$ dataset and various sub-datasets. 
Over all the data, LDA predicts the vanishing order with an accuracy of $95.9\%$ with an explained variance ratio of $0.982$. This means that our linear discriminant is highly accurate at separating the classes and is thus highly informative for the classification of vanishing order.
If we consider a specific subsets, we find some have higher accuracy. 
For example, LDA predictions on just \texttt{G2Q} are $97.1\%$ accurate with an explained variance of $0.997$.

\begin{table}

\begin{center}
\begin{tabular}{c|cccc|l}
Dataset & Training obs. & Validation obs. & Accuracy & Explained Variance & Counts\\
\hline
$\texttt{PRAT}^{\star}$ & \numprint{140924}   & \numprint{35232}  & $0.959$   & $0.982$   &   $0\qquad$ \numprint{53344}\\
&&&&&                                                                $1\qquad$ \numprint{90327} \\
&&&&&                                                                $2\qquad$ \numprint{29648}\\
&&&&&                                                                $3\qquad$ \numprint{2837}\\ \hline\hline
\texttt{BMF} & \numprint{65442}     & \numprint{16361}  & $0.958$   & $0.979$   &   $0\qquad$    \numprint{28280}\\
&&&&&                                                                $1\qquad$ \numprint{44773}\\
&&&&&                                                                $2\qquad$ \numprint{8724}\\
&&&&&                                                                $3\qquad$ \numprint{26}\\
\texttt{ECNF} & \numprint{90791}    & \numprint{22698}  & $0.956$   & $0.983$   &   $0\qquad$ \numprint{42558}\\
&&&&&                                                                $1\qquad$ \numprint{61243}\\
&&&&&                                                                $2\qquad$ \numprint{9661}\\
&&&&&                                                                $3\qquad$ \numprint{27}\\
\texttt{G2Q} & \numprint{50224}     & \numprint{12556}  & $0.971$   & $0.997$   &   $0\qquad$ \numprint{10827}\\
&&&&&                                                                $1\qquad$ \numprint{29155}\\
&&&&&                                                                $2\qquad$ \numprint{19988}\\
&&&&&                                                                $3\qquad$ \numprint{2810} \\
\texttt{HMF} & \numprint{25571}     & \numprint{6393}   & $0.963$   & $0.988$   &    $0\qquad$ \numprint{14443}\\
&&&&&                                                                $1\qquad$ \numprint{16582}\\
&&&&&                                                                $2\qquad$ \numprint{938}\\
&&&&&                                                                $3\qquad$ \numprint{1}\\
\end{tabular}
\caption{LDA results for predicting vanishing order in \texttt{PRAT} and various subsets.}
\label{LDA-table}
\end{center}

\end{table}

\subsection{Neural networks for vanishing order}\label{ss.NN_PRAT}
We may also apply other supervised learning techniques to classify the order of vanishing. 
To that end, we trained two neural networks: one with inputs given by the principal components, 
and the other with inputs given by $v(L)$. 
In both cases, the dataset $\mathcal{D}$ is randomly split into training and test sets exactly as in Section~\ref{s.lda}.

For the network architecture, we tested both feed-forward neural networks (FNN) and 1D convolutional neural networks (CNN). 
We will focus on CNN from this point onward, as it achieved better accuracy than FNN. 
We also experimented with different problem types, namely, classification (outputting the probability of each vanishing order) 
and regression (outputting the estimated value for the vanishing order).
We observe that most of the phenomena appear to be consistent regardless of the network architecture and problem type.

For the hyperparameters of the CNN, we use three 1D convolutional layers, with $16$, $32$ and $64$ channels, 
kernel size $K = 3$ and padding $P = 1$, 
each followed by a $\mathrm{ReLU}$ activation function and Max Pooling layer with $K = 2$ and $P = 1$. 
After a dropout layer $D = 0.5$, the convolutional blocks are followed by two fully connected layers, each with $128$ neurons,
and an output layer whose width is equal to the number of possible orders of vanishing in the dataset.
The networks are trained with cross-entropy loss using Adam optimizer, setting batch size to $3000$ and learning rate to $0.001$.
We use the same set of hyperparameters for both cases.

\begin{figure}[htbp]
    \centering
    \includegraphics[width=0.5\textwidth]{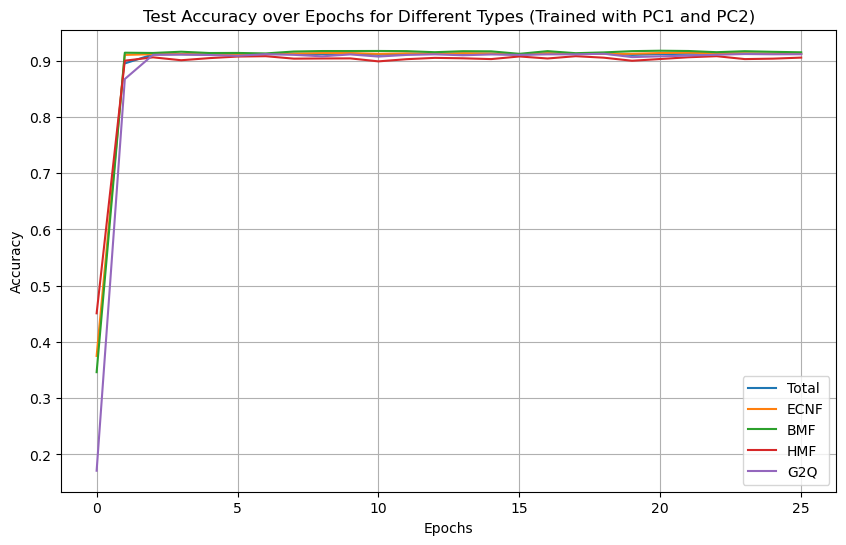}  
    \captionof{figure}{Learning the vanishing order the first and second principal components: percentage accuracy against epoch with CNN. Each colour corresponds to one category of $L$-functions}
    \label{fig:accuracy_PRAT_PC12}
\end{figure}

\begin{figure}[htbp]
    \centering
    \includegraphics[width=0.5\textwidth]{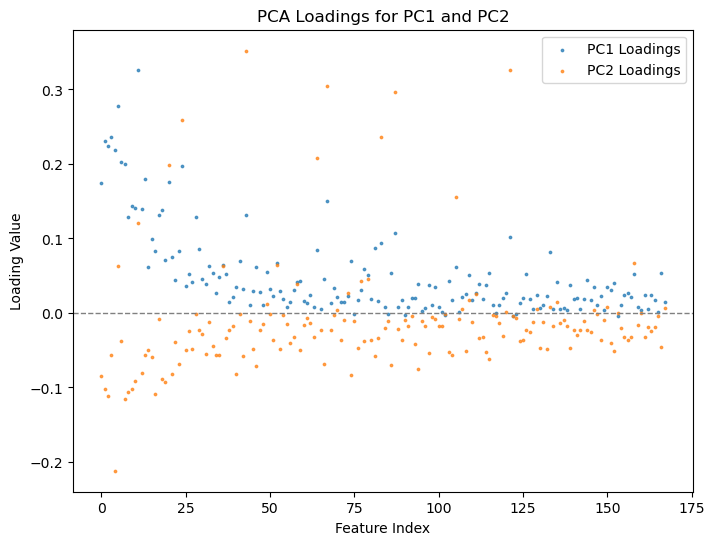} 
    \captionof{figure}{The weights in the first and second principal components used in the training.}
    \label{fig:PC12_loadings}
\end{figure}

When training with principal components,
we observe that a CNN can achieve an $91\%$ overall accuracy 
using only the first two principal components.
The percentage accuracy curves and the final test accuracies
for the named subsets of $\texttt{PRAT}^{\star}$ in this case are presented in 
Figure~\ref{fig:accuracy_PRAT_PC12} and Table~\ref{tab:accuracy_PRAT_PC12}.
In Figure~\ref{fig:PC12_loadings}, we present the weights of each principal component.
We note that similar results were observed in \cite{KV22} using different datasets and features.

\begin{figure}[htbp]
    \centering
    \includegraphics[width=0.5\textwidth]{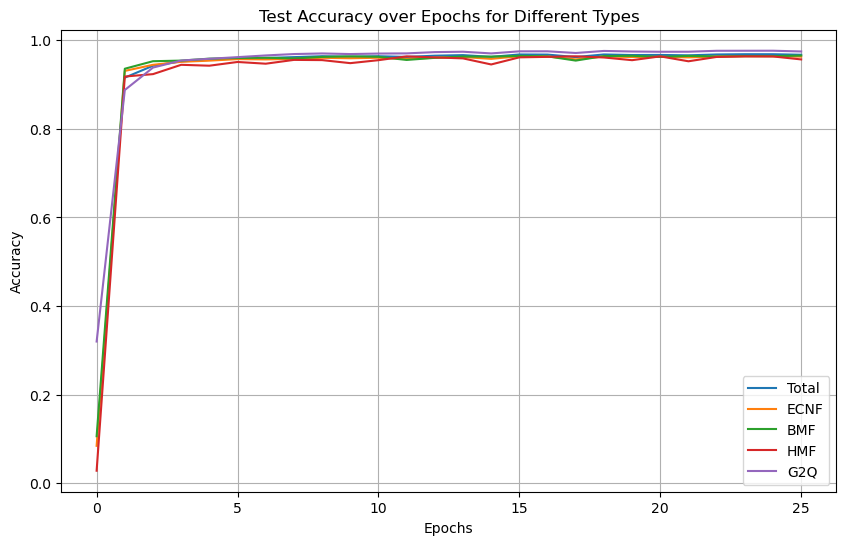} 
    \captionof{figure}{Learning the vanishing order in  $\texttt{PRAT}^{\star}$ from $v(L)$ : percentage accuracy against epoch with CNN. Each colour corresponds to one category of $L$-functions}
    \label{fig:accuracy_PRAT_aps}
\end{figure}

\begin{figure}[htbp]
    \centering
    \begin{subfigure}[b]{0.45\textwidth}
        \centering
        \includegraphics[width=\textwidth]{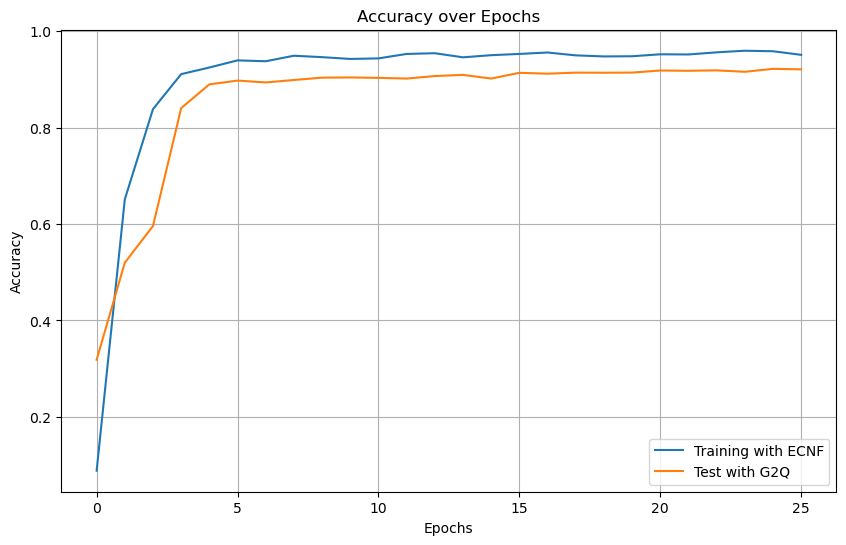} 
    \end{subfigure}
    \hfill
    \begin{subfigure}[b]{0.45\textwidth}
        \centering
        \includegraphics[width=\textwidth]{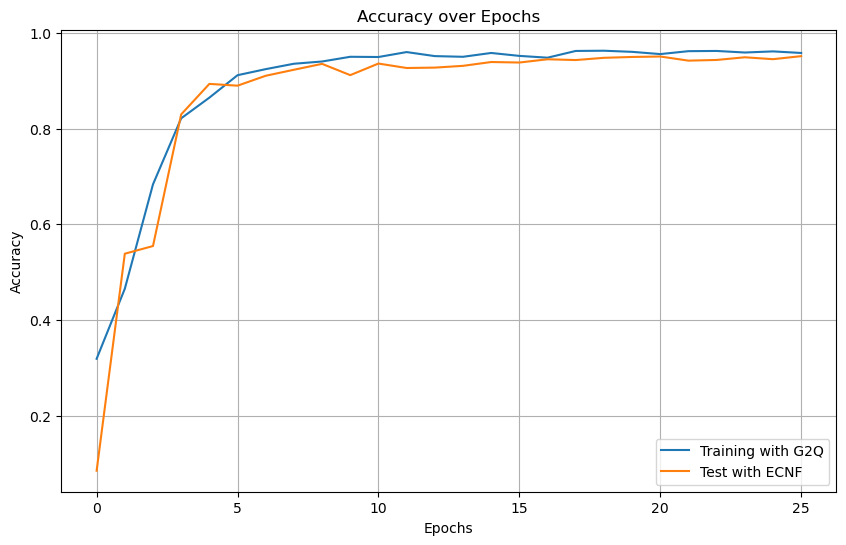}  
    \end{subfigure}
    \caption{Train with \texttt{ECNF} then test with \texttt{G2Q}, and vice versa.}
    \label{fig:transfer_learning}
\end{figure}

When training with $v(L)$, 
the CNN can achieve over 95\% accuracy across all different types of $L$-functions.
The results are shown in Figure~\ref{fig:accuracy_PRAT_aps} and Table~\ref{tab:accuracy_PRAT_aps}.
To show the capability of transfer learning of the neural network,
we first train on only the  \texttt{ECNF} subset  and 
test on \texttt{G2Q}, then reverse the process by training on \texttt{G2Q} and testing on \texttt{ECNF}.
The results are shown in Figure~\ref{fig:transfer_learning},
with both achieving over $90\%$ accuracy on the test set.
This suggests that the features learned from \texttt{ECNF} and \texttt{G2Q} are highly transferable 
and the model generalizes well across $\texttt{PRAT}^{\star}$,
indicating a strong similarity between $L$-functions with the 
same motivic weight and degree.

\begin{table}[htbp]
    \begin{tabular}{c|c}
        sub-dataset & test accuracy\\ \hline
        \texttt{ECNF} & $91.22\%$ \\ 
        \texttt{BMF} & $91.48\%$ \\  
        \texttt{HMF} & $90.54\%$ \\  
        \texttt{G2Q} & $91.13\%$\\   
    \end{tabular}
    \captionof{table}{Test accuracy at the last epoch for each category of $L$-functions in $\texttt{PRAT}^{\star}$, trained with the first and second principal components.}
    \label{tab:accuracy_PRAT_PC12}
\end{table}

\begin{table}[htbp]
    \begin{tabular}{|c|c|}
        \hline
        sub-dataset & test accuracy\\ \hline
        \texttt{ECNF} & $95.37\%$ \\ \hline
        \texttt{BMF} & $95.48\%$ \\ \hline
        \texttt{HMF} & $95.04\%$ \\ \hline
        \texttt{G2Q} & $95.71\%$\\ \hline
    \end{tabular}
    \captionof{table}{Test accuracy at the last epoch for each sub-dataset in $\texttt{PRAT}^{\star}$, trained with $v(L)$.}
    \label{tab:accuracy_PRAT_aps}
\end{table}

\section{Concluding Remarks}\label{s.conc}

The work in this paper gives rise to many open questions about mathematical structure to be found in a data-focused study of rational $L$-functions. 
Predictions using both LDA and neural networks were highly accurate; however, an analysis of the model inaccuracies could give insight into cases where vanishing orders are harder to predict.  

Our work also considered an exploration of transfer learning and augmented training sets, using principal components, for example. We showed that transfer learning is possible, but further exploration is warranted. 
For example, to what extent is the accuracy of our classifiers dependent on the conductor range?
Can we train our model using $L$-functions from a minimal conductor range and then use that model to predict order of vanishing across larger ranges? 
If we train on one data subset, how effective are we at predicting in others? 
Which Dirichlet coefficients are most important in predicting the vanishing order? 
Additionally, how does the prediction of vanishing order trained by principal component results compare to Mestre--Nagao sums?

\appendix
\section{Murmuration Patterns}\label{sect:appendix}

\begin{figure}[htbp!]
    \centering
        \includegraphics[width=.4\textwidth]{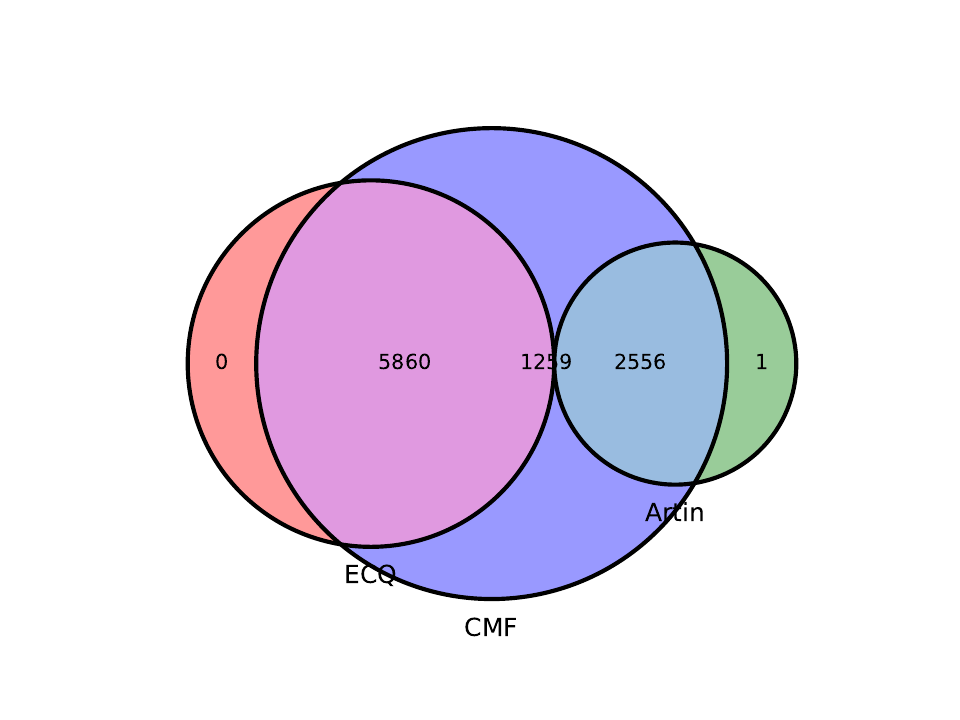}
    \caption{(left) Venn diagram for \texttt{CMF}}
    \label{fig:CMF-venn-murm}
\end{figure}

\begin{figure}[htbp!]
    \centering
    \begin{subfigure}[b]{0.45\textwidth}
        \centering
        \includegraphics[width=\textwidth]{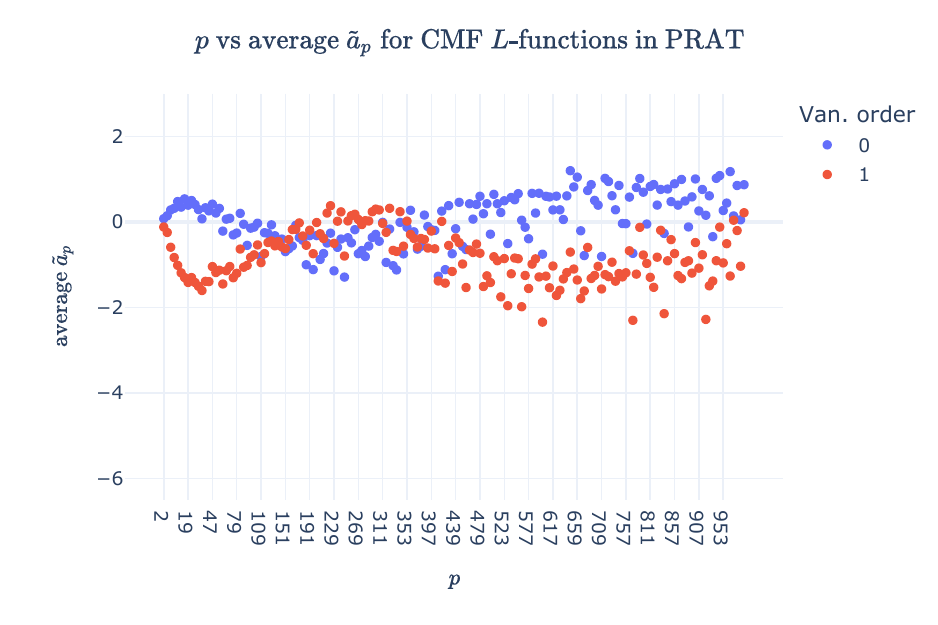}
    \end{subfigure}
    \hfill
    \begin{subfigure}[b]{0.45\textwidth}
        \centering
    \includegraphics[width=\textwidth]{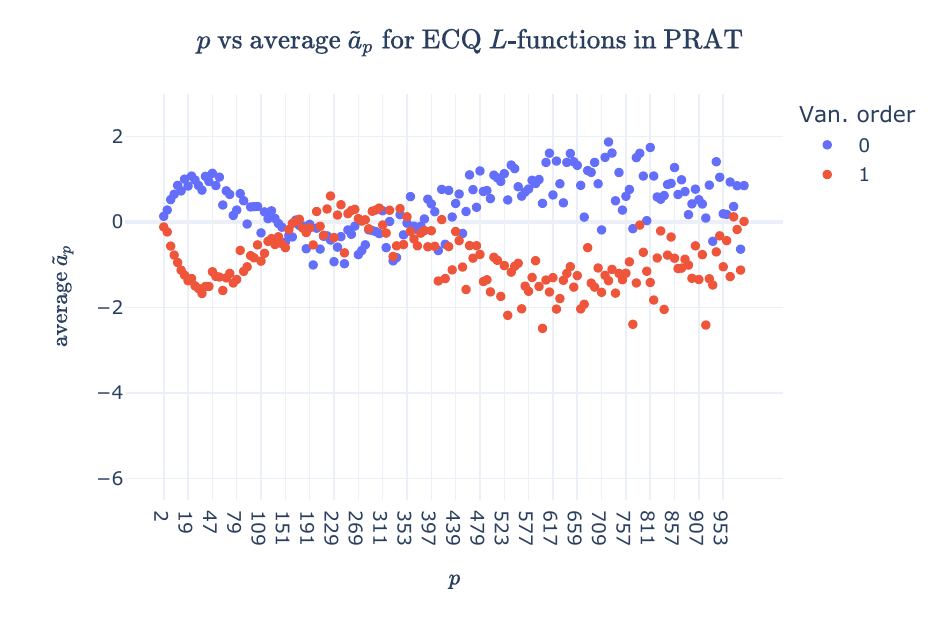}
    \end{subfigure}
    \caption{Average value of $\widetilde{a}_p$ for  (left) \texttt{CMF} and (right) \texttt{ECQ}. 
        }
    \label{fig:murm-ARTIN-ECQ}
\end{figure}

In what follows, we analyze the murmuration patterns for the various sub-datasets in $\texttt{PRAT}_{\leq3}$ that were not contained in  $\texttt{PRAT}^\star$. The murmuration patterns are quite different for these subsets of the data, most likely because the ranges for conductors in the dataset for those $L$-functions are much smaller. 
The $L$-functions not included in $\texttt{PRAT}^\star$ include
\texttt{CMF}, which includes \texttt{ARTIN} and \texttt{ECQ}, and 
\texttt{DIR}. The $L$-functions in \texttt{ARTIN} and \texttt{DIR} all have vanishing order $0$, so we do not include those graphs. Murmurations of Dirichlet characters, grouped by parity (as opposed to vanishing order), were rigorously studied in~\cite{LOPdirichlet}. 
In Figure~\ref{fig:CMF-venn-murm}, we present a Venn diagram for \texttt{CMF}, which includes \texttt{ARTIN} and \texttt{ECQ}.
In Figure~\ref{fig:murm-ARTIN-ECQ}, we plot the average value of $\widetilde{a}_p$ over \texttt{CMF} and \texttt{ECQ}.

Previous work on murmurations has shown that normalizing the $x$-axis by dividing by the conductor is the correct normalisation for plots.
By restricting to root analytic conductor $< 4$ in~\cite{rat}, the conductor range for $L$-functions of degree $2$ is much smaller than the conductor range for $L$-functions of degree $4$.
Thus when we restrict the dataset $\texttt{PRAT}$ to $L$-functions of degree $2$, e.g., \texttt{CMF} and \texttt{ECQ}, we obtain tighter conductor ranges.
This yields murmuration plots with more definition when compared to those plots in Section~\ref{s:murm}.

\bibliographystyle{alpha}
\bibliography{bibfile}

\newcommand{\etalchar}[1]{$^{#1}$}
\begin{thebibliography}{BDKM{\etalchar{+}}13}

\bibitem[AHL{\etalchar{+}}23]{AHLOS}
Malik Amir, Yang-Hui He, Kyu-Hwan Lee, Thomas Oliver, and Eldar Sultanow.
\newblock Machine-learning of the class number of real quadratic fields.
\newblock {\em International Journal of Data Science in the Mathematical
  Sciences}, (1(2)):107--134, 2023.

\bibitem[BBB{\etalchar{+}}21]{best2021computing}
Alex~J. Best, Jonathan Bober, Andrew~R. Booker, Edgar Costa, John~E. Cremona,
  Maarten Derickx, Min Lee, David Lowry-Duda, David Roe, Andrew~V. Sutherland,
  and John Voight.
\newblock Computing classical modular forms.
\newblock In {\em Arithmetic geometry, number theory, and computation}, Simons
  Symp., pages 131--213. Springer, Cham, [2021] \copyright 2021.

\bibitem[BDKM{\etalchar{+}}13]{boberdatabase}
Jonathan Bober, Alyson Deines, Ariah Klages-Mundt, Benjamin LeVeque, R.~Andrew
  Ohana, Ashwath Rabindranath, Paul Sharaba, and William Stein.
\newblock A database of elliptic curves over {$\Bbb Q(\sqrt 5)$}: a first
  report.
\newblock 1:145--166, 2013.

\bibitem[BKN24]{BKN}
Zvonimir Bujanović, Matija Kazalicki, and Lukas Novak.
\newblock Murmurations of mestre-nagao sums.
\newblock \url{http://arxiv.org/abs/2403.17626v1}, 2024.
\newblock {arXiv: math.NT: 2403.17626v1}.

\bibitem[BSS{\etalchar{+}}16]{booker2016database}
Andrew~R. Booker, Jeroen Sijsling, Andrew~V. Sutherland, John Voight, and Dan
  Yasaki.
\newblock A database of genus-2 curves over the rational numbers.
\newblock {\em LMS J. Comput. Math.}, 19:235--254, 2016.

\bibitem[CL25]{rat}
Edgar Costa and The {LMFDB Collaboration}.
\newblock Rational {L}-functions in {LMFDB} with root analytic conductor less
  than 4.
\newblock \url{https://doi.org/10.5281/zenodo.14774042}, 2025.

\bibitem[CN21]{cremona2021q}
J.~E. Cremona and Filip Najman.
\newblock {$\Bbb Q$}-curves over odd degree number fields.
\newblock {\em Res. Number Theory}, 7(4):Paper No. 62, 30, 2021.

\bibitem[Cre81]{cremona1981modular}
J.~E. Cremona.
\newblock {\em Modular symbols}.
\newblock PhD thesis, University of Oxford, 1981.

\bibitem[Cre84]{cremona1984hyperbolic}
Jonh~E. Cremona.
\newblock Hyperbolic tessellations, modular symbols, and elliptic curves over
  complex quadratic fields.
\newblock {\em Compositio Math.}, 51(3):275--324, 1984.

\bibitem[Cre92]{cremona1997algorithms}
J.~E. Cremona.
\newblock {\em Algorithms for modular elliptic curves}.
\newblock Cambridge University Press, Cambridge, 1992.

\bibitem[DV21]{donnelly2021database}
Steve Donnelly and John Voight.
\newblock A database of {H}ilbert modular forms.
\newblock In {\em Arithmetic geometry, number theory, and computation}, Simons
  Symp., pages 365--373. Springer, Cham, [2021] \copyright 2021.

\bibitem[HLO22a]{HLOb}
Yang-Hui He, Kyu-Hwan Lee, and Thomas Oliver.
\newblock Machine-learning number fields.
\newblock {\em Math. Comput. Geom. Data}, 2(1):49--66, 2022.

\bibitem[HLO22b]{HLOa}
Yang-Hui He, Kyu-Hwan Lee, and Thomas Oliver.
\newblock Machine-learning the {S}ato-{T}ate conjecture.
\newblock {\em J. Symbolic Comput.}, 111:61--72, 2022.

\bibitem[HLO23]{HLOc}
Yang-Hui He, Kyu-Hwan Lee, and Thomas Oliver.
\newblock Machine learning invariants of arithmetic curves.
\newblock {\em J. Symbolic Comput.}, 115:478--491, 2023.

\bibitem[HLOP24]{HLOP}
Yang-Hui He, Kyu-Hwan Lee, Thomas Oliver, and Alexey Pozdnyakov.
\newblock Murmurations of elliptic curves.
\newblock {\em Experimental Mathematics}, pages 1--13, 2024.

\bibitem[JR17]{jones2017artin}
John~W. Jones and David~P. Roberts.
\newblock Artin {$L$}-functions of small conductor.
\newblock {\em Res. Number Theory}, 3:Paper No. 16, 33, 2017.

\bibitem[KV23]{KV22}
Matija Kazalicki and Domagoj Vlah.
\newblock Ranks of elliptic curves and deep neural networks.
\newblock {\em Res. Number Theory}, 9(3):Paper No. 53, 21, 2023.

\bibitem[LGS{\etalchar{+}}14]{lex2014upset}
Alexander Lex, Nils Gehlenborg, Hendrik Strobelt, Romain Vuillemot, and
  Hanspeter Pfister.
\newblock {UpSet}: visualization of intersecting sets.
\newblock {\em IEEE transactions on visualization and computer graphics},
  20(12):1983--1992, 2014.

\bibitem[{LMF}24]{lmfdb}
The {LMFDB Collaboration}.
\newblock The {L}-functions and modular forms database.
\newblock \url{https://www.lmfdb.org}, 2024.
\newblock [Online; accessed 27 November 2024].

\bibitem[LOP25]{LOPdirichlet}
Kyu-Hwan Lee, Thomas Oliver, and Alexey Pozdnyakov.
\newblock Murmurations of {D}irichlet {C}haracters.
\newblock {\em Int. Math. Res. Not. IMRN}, (1):rnae277, 2025.

\bibitem[Oli24]{O24}
Thomas Oliver.
\newblock Machine learning for number theory: unsupervised learning with
  {$L$}-functions.
\newblock In {\em Mathematical software---{ICMS} 2024}, volume 14749 of {\em
  Lecture Notes in Comput. Sci.}, pages 196--203. Springer, Cham, [2024]
  \copyright 2024.

\bibitem[Poz24a]{Poz}
Alexey Pozdnyakov.
\newblock Predicting root numbers with neural networks.
\newblock \url{http://arxiv.org/abs/2403.14631v1}, 2024.
\newblock {arXiv:math.NT:2403.14631v1}.

\bibitem[Poz24b]{alexeypozdnyakov2024predicting}
Alexey Pozdnyakov.
\newblock Predicting root numbers with neural networks.
\newblock \url{http://arxiv.org/abs/2403.14631v1}, 2024.
\newblock {arXiv:math.NT:2403.14631v1}.

\bibitem[Sar23]{sarnakLetter}
Peter Sarnak.
\newblock Letter to sutherland and zubrilina.
\newblock \url{https://publications.ias.edu/sarnak/paper/2726}, 2023.

\bibitem[Sel92]{S92}
Atle Selberg.
\newblock Old and new conjectures and results about a class of {D}irichlet
  series.
\newblock In {\em Proceedings of the {A}malfi {C}onference on {A}nalytic
  {N}umber {T}heory ({M}aiori, 1989)}, pages 367--385. Univ. Salerno, Salerno,
  1992.

\bibitem[Sut22]{drewLetter}
Andrew Sutherland.
\newblock Letter from sutherland to rubinstein and sarnak.
\newblock \url{https://math.mit.edu/~drew/RubinsteinSarnakLetter.pdf}, August
  2022.

\bibitem[Whi90]{whitley1990modular}
Elise Whitley.
\newblock {\em Modular forms and elliptic curves over imaginary quadratic
  number fields}.
\newblock PhD thesis, University of Exeter, 1990.

\end{thebibliography}

\end{document}